\documentclass{article}
\usepackage[all]{xy} 
\usepackage{latexsym,amsmath,amssymb}
\usepackage{fancyhdr}
\usepackage[english]{babel}
\usepackage{graphicx}
\usepackage{verbatim}
\title{\bf Filter convergence and decompositions for vector lattice-valued measures}
\author{\bf D. Candeloro --A. R. Sambucini
\footnote{ \noindent Supported by GNAMPA of CNR and University of Perugia \newline
corresponding author: Domenico Candeloro \newline
Authors' Address: Dipartimento di Matematica e Informatica, via Vanvitelli,1
 I-06123 Perugia (Italy)
e-mail: candelor@dmi.unipg.it;anna.sambucini@unipg.it}}
\date{\Large   28 gennaio 2014}
\begin{document}
\maketitle \pagestyle{myheadings} \markboth{\centerline{\small \sl
Convergence and decompositions}}{\centerline{\small \sl Convergence and decompositions}}
\newtheorem{teorema}{Theorem}[section]
\newtheorem{tteorema}[teorema]{\Large Theorem}
\newtheorem{lemma}[teorema]{Lemma}
\newtheorem{proposizione}[teorema]{Proposition}
\newtheorem{corollario}[teorema]{Corollary}
\newtheorem{osservazione}[teorema]{Remark}
\newtheorem{remark}[teorema]{Remark}
\newtheorem{remarks}[teorema]{Remarks}
\newtheorem{definizione}[teorema]{Definition}
\newtheorem{condizione}[teorema]{Condition}
\newtheorem{notazioni}[teorema]{Notations}
\newtheorem{examples}[teorema]{Examples}
\newtheorem{example}[teorema]{Example}
\newtheorem{definizioni}[teorema]{Definitions}
\newtheorem{definitions}[teorema]{Definitions}
\newtheorem{assunzioni}[teorema]{Assumptions}
\newcommand{\erre}{\mbox{$\mathbb{R}$}}
\newcommand{\enne}{\mbox{$\mathbb{N}$}}
\newcommand{\bvi}{\bigvee_{i=1}^{\infty} a_{i, \varphi(i)}}
\newcommand{\bvib}{\bigvee_{i=1}^{\infty} b_{i, \varphi(i)}}
\newcommand{\bvic}{\bigvee_{i=1}^{\infty} c_{i, \varphi(i)}}
\newcommand{\bvid}{\bigvee_{i=1}^{\infty} d_{i, \varphi(i)}}
\newcommand{\bvie}{\bigvee_{i=1}^{\infty} e_{i, \varphi(i)}}
\newcommand{\bviA}{\bigvee_{i=1}^{\infty} A_{i, \varphi(i)}}
\newcommand{\bviB}{\bigvee_{i=1}^{\infty} B_{i, \varphi(i)}}
\newcommand{\Fs}{\mbox{$\mathcal F$}}
\newcommand{\Hs}{\mbox{$\mathcal H$}}
\newcommand{\Ps}{\mbox{$\mathcal P$}}
\newcommand{\Gs}{\mbox{$\mathcal G$}}
\newcommand{\Is}{\mbox{$\mathcal I$}}
\newcommand{\As}{\mbox{$\mathcal A$}}
\newcommand{\Bs}{\mbox{$\mathcal B$}}
\newcommand{\ttilde}{~}
%_________________________________________________________________________
\newcommand{\Ffreccia}{\mbox{$ \stackrel{\,\,\,\mathcal F\,\,\,}{\rightarrow}$}}
\newcommand{\uFfreccia}{\mbox{$ \stackrel{\,\,\,u_{\mathcal F}\,\,\,}{\rightarrow}$}}
\newcommand{\rFfreccia}{\mbox{$ \stackrel{\,\,\,r_{\mathcal F}\,\,\,}{\rightarrow}$}}
\newcommand{\oFfreccia}{\mbox{$ \stackrel{\,\,\,o_{\mathcal F}\,\,\,}{\rightarrow}$}}
\newcommand{\DFfreccia}{\mbox{$ \stackrel{\,\,\,D_{\mathcal F}\,\,\,}{\rightarrow}$}}
\newcommand{\sFfreccia}{\mbox{$ \stackrel{\,\,\, \sigma_{\mathcal F}\,\,\,}{\rightarrow}$}}
%_________________________________________________________________________
\centerline{\em Dedicated to Prof. Benedetto Bongiorno}
\bigskip
\begin{abstract} 
Filter convergence of vector lattice-valued measures is considered, in order to deduce theorems of convergence for their decompositions. First the $\sigma$-additive case is studied, without particular assumptions on the filter; later the finitely additive case is faced, first assuming uniform $s$-boundedness (without restrictions on the filter), then relaxing this condition but imposing stronger properties on the filter. In order to obtain the last results, a Schur-type convergence theorem, obtained in \cite{bdpschurfilters}, is used.
\end{abstract}
\mbox{~} \vskip.1cm
\noindent
{\bf  2010 AMS Mathematics Subject Classification}: {\rm 28B15, 28B05, 06A06, 54F05.}\\
{\bf Key words}: \rm  Vector lattices, filter convergence, Lebesgue decomposition, Sobczyk-Hammer decomposition, Yosida-Hewett decomposition.
\normalsize

\bigskip

\section{Introduction}
The wide literature originated by the classic Nikody'm and Vitali-Hahn-Saks theorems has produced interesting results in many abstract settings (see for instance \cite{mimmo, 4, DSNP, cdl}). In this framework recently also the so-called {\em filter convergence} has been considered, though in this case it is quite difficult to obtain general positive results. 
Filter convergence (or, equivalently, ideal convergence) has been extensively studied in the literature, since the papers  \cite{KSW} and \cite{NR}.  Besides real-valued measures, also topological group-valued and vector lattice-valued measures have been considered recently: we recall here \cite{Dallas,bd2} for applications
in abstract integration theory and approximation theory,
and \cite{cascales, bdpschurfilters,bd3} for convergence results related to the Schur theorem.
Another interesting problem related to measures consists in the possibility of decomposing them into {\em simpler} parts (Lebesgue decomposition, and so on). Many questions of this type have been considered recently, also in very abstract settings: we just mention \cite{DAG, DSN1, DSN2, ABVW1, AV2, W1, bongio}.

 Our research considers filter convergence of vector lattice-valued measures (both finitely and $\sigma$-additive ones) in order to deduce convergence of their decompositions, in the sense of Lebesgue,  Sobczyk-Hammer, or Yosida-Hewett. When the vector lattice is super-Dedekind complete and weakly-$\sigma$-distributive, in the $\sigma$-additive case convergence of the decompositions is almost trivial, thanks to well-known results concerning the decompositions themselves (see \cite{prace, candeloro}); while in the finitely additive setting a kind of uniform $s$-boundedness is needed, which makes the task considerably harder. The main results here obtained are Theorem \ref{main} and Theorem \ref{finale}, together with the Corollary \ref{corofinale}: thanks to a Schur-type theorem obtained in \cite{ bdpschurfilters}, we are able to prove that, for a suitable kind of filters, it is possible to obtain uniform $s$-boundedness from a weaker property, called {\em ideal uniform $s$-boundedness}, which is always satisfied when the filter is the natural one of cofinite sets in $\enne$, usually denoted by $\Fs_{cofin}$.
Finally an alternative version of  Theorem \ref{finale} is given in Theorem \ref{nuovoschur},
relaxing some conditions on the domain of the measures but strenghtening conditions on the filter.
Part of the results here obtained have been presented at the Mini-Symposium on Real Analysis, Measure Theory and Integration, held in Palermo, May 13, 2013, and dedicated to Prof. Benedetto Bongiorno in occasion of his $70^{th}$ birthday.
%________________________________________________________________________________________
%_______________________________________________________________________________________-
\section{Vector lattice-valued measures}

Here we shall assume that
${\bf X}$ is a vector lattice, i.e. $({\bf X},+,\cdot,\leq)$ is a linear space endowed with a compatible order, stable w.r.t. finite suprema and infima.

We assume also that ${\bf X}$ is {\em super-Dedekind complete}, i.e.  
  any subset $A\subset {\bf X}$ bounded from above has {\em supremum} in ${\bf X}$, and
 $\vee(A)=\vee(A_0)$
for some {\em countable}  subset $A_0\subset A$. 

\begin{definizioni}\rm
We call {\em $(o)$-sequence}
 any decreasing sequence $(r_n)_n$ in ${\bf X}^+$ such that $\wedge r_n=0$.
\\
 A sequence  $(x_k)_k$ is { \em $(o)$-convergent} to $x\in {\bf X}$ if the following is an $(o)$-sequence:
$$r_n:=\bigvee_{k\geq n}|x_k-x|.$$
We shall call {\em $D$-sequence} 
any bounded double sequence $(a_{i,j})_{i,j}$ in ${\bf X}^+$  such that
$(a_{i,j})_j$ is an $(o)$-sequence
  for every fixed $i$.
\\
 Given a $D$-sequence $(a_{i,j})_{i,j}$, and any mapping $\varphi:\enne \to \enne$, the element
$$v:=\bigvee_i a_{i,\varphi(i)}$$ is said to be a
 {\em domination} for $(a_{i,j})_{i,j}$. 
\\
 A sequence $(r_n)_n$ in ${\bf X}$ is { \em $(D)$-convergent} to $r$ if there exists a $D$-sequence $(a_{i,j})$ whose every domination is definitely larger than $|r_n-r|$.\\
A vector lattice ${\bf X}$ is {\em weakly $\sigma$-distributive} if, for every $D$-sequence $(a_{i,j})_{i,j}$, the set of its dominations has infimum 0.
\end{definizioni}

From now on  we shall always assume that 
${\bf X}$ is a {\bf super-Dedekind complete} and {\bf weakly $\sigma$-distributive} vector lattice.\\

Usually we shall denote by $\Omega$ any abstract space, by $\Hs$ any algebra of subsets of $\Omega$, and by $\As$ any $\sigma$-algebra in $\Omega$. 

Now we shall consider finitely additive measures $m:\Hs\to {\bf X}$. The following definitions have been already introduced in the paper \cite{prace}, though in terms of $D$-sequences: in our setting these notions are equivalently stated as follows, by means of $(o)$-sequences.

\begin{definizioni}{\label{basiche}}\rm
We say that $m$ is {\em $s$-bounded} if there exists an $(o)$-sequence $(p_n)$ in ${\bf X}$ such that, for every pairwise disjoint sequence $(F_k)_k$ from $\Hs$ and every integer $n$, it is possible to find an index $k(n)$ satisfying
\begin{eqnarray}
\bigvee_{k\geq k(n)}|m(F_{k})|\leq p_n.
\end{eqnarray}
If this is the case, we also say that the $(o)$-sequence $(p_n)$ {\em regulates} $s$-boundedness of $m$.\\
\noindent
An $s$-bounded finitely additive measure $m$ is also bounded, (see
 {\rm \cite[Theorem 5.3]{4}}) and in this case it is possible to define its {\em variations}:
$$v^+(m) (H)=\sup_{A\in {\mathcal H}}m(A\cap H),\ \ v^-(m)(H)=-\inf_{A\in  {\mathcal H}}m(A\cap H),\ \ v(m)=v^++v^-.$$
It turns out that, if $m$ is bounded, then $m$ is $s$-bounded if and only if its {\em total variation} $v(m)$ is.
\end{definizioni}

\begin{definizione} \rm
We say that $m$ is {\em $\sigma$-additive} if there exists an $(o)$-sequence $(p_n)$ in ${\bf X}$ such that, 
for every decreasing sequence $(F_k)_k$ from $\Hs$ with empty
 intersection, and every integer $n$ it is possible to find an index $k(n)$ satisfying
\begin{eqnarray}
\bigvee_{A\in  {\mathcal H}}|m(A\cap F_{k(n)})|\leq p_n.
\end{eqnarray}
Also in this case, we say that the $(o)$-sequence $(p_n)$ {\em regulates} $\sigma$-additivity of $m$.
\end{definizione}

In this regard we recall an elegant technique that is often used in order to  study finitely additive measures by means of {\em isomorphic} $\sigma$-additive ones, i.e. the Stone Isomorphism method: details can be found in \cite{4,panoramica}. As it is well-known, \\

\begin{description}\label{iso}
\item[($H_{\mbox{\tiny Stone}}$)] given $\Omega$ and $\Hs$, there exist a compact, totally disconnected Hausdorff space $S$, and an algebraic isomorphism $\psi:\Hs\to \Sigma$ where $\Sigma$ is the algebra of all clopen subsets of $S$. As usual, the space $S$ will be called the {\em Stone space} associated with $\Hs$.\\
\end{description}

We state here one of its formulations, which is valid when ${\bf X}$ is super-Dedekind complete and weakly-$\sigma$-distributive.

\begin{teorema}{\label{stone}}
Let $m:\Hs\to {\bf X}$ be any ${\bf X}$-valued $s$-bounded finitely additive measure.
\begin{description}
\item[\ref{stone}.1] The measure $m\cdot \psi^{-1}:\Sigma\to {\bf X}$ has a unique $\sigma$-additive extension $\widetilde{m}$ to the $\sigma$-algebra $\sigma(\Sigma)$ generated by $\Sigma$, and
\item[\ref{stone}.2] there exists an $(o)$-sequence $(p_n)_n$ in ${\bf X}$ such that, for every element $A\in \sigma(\Sigma)$ and every integer $n$ there corresponds an element $F_n\in \Hs$ satisfying
\begin{eqnarray*}
\sup_{F\in {\mathcal F}}|\widetilde{m}(\psi(F)\cap(A  \, \Delta\, \psi(F_n))|\leq p_n.
\end{eqnarray*}
\end{description}
\end{teorema}

\begin{osservazione}\rm
The measure $\widetilde{m}$  is 
called the 
{\em Stone extension} of $m$ (though it {\em lives} in a different space), and the last property, {\bf \ref{stone}.2}, will be mentioned as the {\em density property } of $\Sigma$  (or also of $\Hs$) in $\sigma(\Sigma)$ with respect to $\widetilde{m}$.

\end{osservazione}

Another useful technique allows to {\em transform} finite additivity into $\sigma$-additivity. 
\begin{teorema}{\label{restrizioni}}{\rm \cite[Theorem 5.3]{prace}}
Let $m:\As\to {\bf X}$ be any  $s$-bounded finitely additive measure, defined on a $\sigma$-algebra $\As$, and let $(H_n)_n$ be any pairwise disjoint family from $\As$. Then there exists a sub-sequence $(H_{n_k})_k$ such that $m$ is $\sigma$-additive in the $\sigma$-algebra generated by the sets  $H_{n_k}$. Moreover, if $(r_n)_n$ denotes any $(o)$-sequence regulating $s$-boundedness of  $v(m)$, then the same $(o)$-sequence regulates $\sigma$-additivity along the subsequence $H_{n_k}$. 
\end{teorema}

\begin{osservazione}{\label{processodiagonale}}\rm
By means of a diagonal argument, it is possible to prove that, in case of a {\em sequence} of $s$-bounded finitely additive measures $(m_j)$, whenever $(H_n)_n$ is any pairwise disjoint family from $\As$, a unique subsequence $(H_{n_k})_k$ can be found, such that all the measures $(m_j)$ are $\sigma$-additive along the  $\sigma$-algebra generated by the sets  $H_{n_k}$. Also in this assertion, we can observe that, in case $(r_n)_n$ is an $(o)$-sequence regulating $s$-boundedness of {\em all} the measures $v(m_j)$, then the same $(o)$-sequence regulates $\sigma$-additivity of all measures along the subsequence $(H_{n_k})_k$.
\end{osservazione}

%________________________________________________________________________________________
%_______________________________________________________________________________________-
\section{Decompositions}

Assume  that, together with $m:\Hs\to {\bf X}$, also a positive real-valued finitely additive measure $\nu:\Hs\to \erre$ is fixed. Then we say that $m$ is {\em absolutely continuous} with respect to $\nu$ (and write $m\ll \nu$) when the following setting defines an $(o)$-sequence in ${\bf X}$: 
$$p_n:=\sup\{|m(A)|:A\in \Hs, \nu(A)\leq \frac{1}{n}\},\ \ \ \ \ \ n\in \enne.$$

In case $m$ and $\nu$ are $\sigma$-additive and $\Hs$ is a $\sigma$-algebra, then an equivalent definition (see \cite{4,panoramica}) of absolute continuity is as follows (for $A\in \Hs$):
$$\nu(A)=0\Rightarrow m(A)=0.$$

Finally, assuming that $m:\Hs\to {\bf X}$ and $\nu:\Hs\to \erre^+_0$ are two finitely additive measures, we say that $m$ and $\nu$ are mutually {\em singular} (and we write $\nu\perp m$) if there exist a  sequence $(A_k)_k$ in $\Hs$ and an $(o)$-sequence $(q_k)_k$ in ${\bf X}$ such that $\lim_k\nu(A_k)=0$ and, for every $k$ 
$$\sup\{|m(E\setminus A_k)|:E\in \Hs\}\leq q_k.$$

We also remark that, in case $m$ and $\nu$ are $\sigma$-additive on a $\sigma$-algebra $\Hs$, singularity is equivalent to the existence of a set $F\in \Hs$ such that $m(A\cap F)=\nu(F^c)=0$ for every $A\in \Hs$.

A well-known decomposition theorem for measures can be found in \cite{candeloro}. We first remark that, given any finitely additive measure $m$ on some algebra $\Hs$ and a set $F\in \Hs$, the notation $m_{|F}$ stands for the measure $A\mapsto m(A\cap F)$, for all $A\in \Hs$.

\begin{teorema}{\label{troncate}}{\rm (\cite[Theorem 5.2]{candeloro})}
Let $m:\As\to {\bf X}$ and $\nu:\As\to \erre^+_0$ be two $\sigma$-additive measures on a $\sigma$-algebra $\As$. Then there exists an element $E\in \As$ such that 
$$m_{|E}\ll \nu,\ \ \ m_{|E^c}\perp \nu.$$
\end{teorema}
As usual, the measures $m_{|E}$ and $m_{|E^c}$ are called the {\em absolutely continuous} and {\em singular} parts of $m$ w.r.t. $\nu$, and denoted by $m^{<}$ and $m^{\perp}$ respectively.\\

This kind of decomposition has been extended to the finitely additive case too, by using the Stone extension method, which is valid for $s$-bounded finitely additive ${\bf X}$-valued measures (see {\rm (\cite[Theorems 2.12 and 4.8]{panoramica})}).   
Similarly one can proceed to obtain other kinds of decompositions, e.g Yosida-Hewett and Sobczyk-Hammer decompositions.

\begin{definizioni}{\label{cont}}\rm
Let $m:\Hs\to {\bf X}^+_0$ be any finitely additive measure. We say that $m$ is {\em continuous} if there exists an $(o)$-sequence $(p_n)_n$ in ${\bf X}$ and a sequence $(\pi_n)_n$ of finite partitions of $\Omega$, $\pi_n=\{J_1,...,J_{k_n}\}$ such that for each $n$ we have
$$\sup_{i=1...,k_n}m(J_i)\leq p_n.$$
(Here {\em partition} means that all sets $J_i$ are in $\Hs$, are pairwise disjoint and their union is $\Omega$). 

We also say that $m$ is {\em atomic} if there exist no continuous finitely additive measure $\mu:\Hs\to {\bf X}^+_0$ such that $\mu\leq m$.

In case $m$ is not a positive measure, but is $s$-bounded, then it will be said to be {\em continuous} if its total variation $v(m)$ is.
\end{definizioni}
Also for this notion we have a decomposition theorem, see \cite{panoramica}. We give here the formulation in the $\sigma$-additive case, for the sake of simplicity. 

\begin{teorema}{\label{sobham1}}{\rm (\cite[Theorem 5.6]{panoramica})}
Let us assume that $m:\As \to {\bf X}^+_0$ is any $\sigma$-additive measure on a $\sigma$-algebra $\As$. Then there exists an element 
$F\in \As$ such that $m_{|F}$ is continuous and $m_{|F^c}$ is atomic.
\end{teorema}

As usual, the measures $m_{|F}$ and $m_{|F^c}$ are called the {\em continuous} and {\em atomic} parts of $m$, and denoted by $m^{s}$ and $m^{a}$ respectively.

Another kind of decomposition is considered often, i.e. the so-called Yosida-Hewett one, for finitely additive measures on a $\sigma$-algebra. 

\begin{definizione}{\label{yosidahewett}}\rm
Given a non-negative finitely additive measure $m:\As\to {\bf X}$, defined on a $\sigma$-algebra $\As$, we say that $m$ is {\em purely finitely additive} if $0$ is the only countably additive measure $\mu:\As\to {\bf X}^+_0$,  such that $\mu\leq m$.
\end{definizione}

The following result holds true, see \cite{panoramica}.

\begin{teorema}{\label{yosidhewe}}{\rm (\cite[Theorem 6.3]{panoramica})}
Given a  finitely additive measure $m:\As\to {\bf X}_0^+$, defined on a $\sigma$-algebra $\As$, there exists a unique purely finitely additive measure $m_0:\As\to {\bf X}^+_0$, such that $m_0\leq m$, and such that $m-m_0$ is countably additive.
\end{teorema}

The measure $m_0$ can be described as follows:  with the notations in 
{\bf ($H_{\mbox{\tiny Stone}}$)}
and Theorem \ref{stone}, there exists a set $U$ in the $\sigma$-algebra $\sigma(\Sigma)$ such that $m_{:U}$ is purely finitely additive and $m_{:U^c}$ is $\sigma$-additive, where $m_{:U}(A):=\widetilde{m}(\psi(A)\cap U)$, for $A\in \As$.
%________________________________________________________________________________________
%_______________________________________________________________________________________- 
\section{Filter convergence} 
We briefly recall some definitions and properties  related to the so-called {\em filter convergence}: further details can be found in \cite{bdpschurfilters,Dallas}, though in terms of $D$-sequences rather than $(o)$-sequences (which in our context is perfectly equivalent).

Let $Z$ be any fixed  set. 
\begin{definizione}\label{filter} \rm
A family  $\mathcal{F}$ of subsets of $Z$
is called  a { \textit{filter}} in $Z$ iff 
\begin{description}
\item[\ref{filter}.a] $\emptyset  \not \in {\mathcal F}$,
\item[\ref{filter}.b] $A \cap B \in {\mathcal F}$ whenever $A,B \in \mathcal{F}$ and
\item[\ref{filter}.c] $A \in \mathcal{F}, B \supset A \Rightarrow B\in \mathcal{F}.$
\end{description}
\end{definizione}
Given any filter $\mathcal{F}$ of subsets of $Z$, the {\em dual ideal} of $\mathcal{F}$ is the family of all {\em complements} of elements from $\mathcal{F}$. Usually the dual ideal will be denoted by $\mathcal{I}_{ \mathcal{F}}$.\\ 
If the dual ideal $\mathcal{I}_{ \mathcal{F}}$ contains all finite subsets of $Z$, we say that  $\mathcal{F}$ is a {\em free} filter. 
Given any filter $\mathcal{F}$  in $Z$, we say that a subset $H\subset Z$ is {\em stationary} if it does not belong to the dual ideal $\mathcal{I}_{ \mathcal{F}}$, i.e. if and only if $H\cap F\neq \emptyset$ for all $F\in \mathcal{F}$.\\

Let ${\mathcal F}$ be any filter in $\enne$.
 \begin{definizioni}\rm
 For every infinite set $I\subset \enne$ a {\em block} of $I$ is any partition $\{D_k,k\in \enne\}$ of $I$, obtained with {\em finite} sets $D_k$ in $\enne$. \\
 The filter $\Fs$ is said to be {\em block-respecting} if, for every stationary set $H$ and every block $\{D_k:k\in \enne\}$ of $H$ there exists a stationary set $J\subset H$ such that 
$\left\vert{J\cap D_k}\right\vert\leq 1$ for all $k$.
\\
The filter $\mathcal{F}$ is said to be {\em diagonal} if for every
sequence $(A_n)_n$ in $\mathcal{F}$ and every stationary set $I\subset \enne$, there exists a stationary set $J\subset I$ such that the set $J\setminus A_n$ is finite for all $n\in \enne$.
 \end{definizioni}
 
From now on, we assume that our involved filters are free.

\medskip

\noindent We want to stress the fact that there exist filters enjoying both properties of {\em block respecting} and {\em diagonality}, other than $\Fs_{cofin}$. 
For example, any {\em countably generated} filter has these properties (see {\rm \cite[Sec. 4]{cascales}}): for further reference, we give here an equivalent definition, based on the dual ideal.
\begin{definizione}\label{4.3}\rm
Given an ideal $\Is$  in $\Ps(\enne)$, we say that $\Is$ is {\em countably generated} if there exists a partition $(A_k)_k$ of $\enne$ into pairwise disjoint nonempty subsets, such that $\Is$ is precisely the set of all finite unions of subsets of the sets $A_k$.
\end{definizione}
Of course, a filter $\Fs$ is countably generated if and only if its dual ideal is.\\

Later, we shall need the following 

\begin{proposizione}\label{superdiagonal}
Let $\Fs$ be a countably generated filter. Then, for every stationary set $J\in \Fs$ there exists a stationary subset $J'\subset J$, such that every infinite subset of $J'$ is stationary.

\end{proposizione}
{\bf Proof.}\ Let $(A_n)_n$ denote a partition as in Definition \ref{4.3}, for the dual ideal of $\Fs$. Since $J$ is stationary, the set $J\cap A_n$ is nonempty for infinitely many values of $n$. Picking one point in each of the nonempty intersections, we obtain an infinite set $J'\subset J$. Since $\left\vert{(J'\cap A_n)}\right\vert\leq 1$ for each $n$, 
$J'$ is not in the dual ideal, i.e. is stationary. For the same reason, every infinite subset of $J'$ is stationary too.\ \  $\Box$

\begin{definizione}\label{filterconv} \rm

A sequence $(x_k)_{k \in \enne}$ in ${\bf X}$ 
{ \em{$(o_{\mathcal F})$-converges to $x
\in {\bf X}$} }   ($x_k \oFfreccia x$) iff
 there exists 
an $(o)$-sequence $(\sigma_p)_p $ in ${\bf X}$  such
that  the set 
$$\{ k \in \enne: |x_k - x| \leq \sigma_p  \}$$
is an element of ${\mathcal F}$ for each $p \in \mathbb{N}$.
\end{definizione}

%________________________________________________________________________________________
%_______________________________________________________________________________________-
\subsection{ Convergence for decompositions in the $\sigma$-additive case}

In this section we shall prove that, independent of the filter involved, filter pointwise convergence of the sequence $(m_n)$ essentially implies the same for the corresponding sequences of decompositions.
As before, we shall assume that there exists a fixed positive $\sigma$-additive measure $\nu$ defined on the $\sigma$-algebra $\As$.

More precisely we have the following
\begin{teorema}{\label{teokyber}}
If $(m_n)_n$ is any sequence of bounded ${\bf X}$-valued $\sigma$-additive measures defined on a measure space $(\Omega,\As,\nu)$ 
such that
 the sequence $(m_n)_n$ is pointwise $(o_{\mathcal F})$-convergent to a { $\sigma$-additive} measure $\mu$.
Then the sequences 
$$(m_n^{<})_n,\ (m_n^{\perp})_n, \ (m_n^a)_n,\ \ (m_n^s)_n$$ converge in the same way to
$m^{<},\ m^{\perp}, \ m^a,\ \ m^s$ 
 respectively.
\end{teorema}

{\bf Proof.}\  The proof  is essentially based upon theorems \ref{troncate} and \ref{sobham1} : given a sequence $(m_n)_n$ of ${\bf X}$-valued $\sigma$-additive measures on a $\sigma$-algebra $\As$, and a real-valued positive $\sigma$-additive measure $\nu$ on  $\As$, then there exists a {\em unique} set $U\in \As$ such that 
$$ \ \ \ {m_n}_{|U}\perp \nu, \ \ { m_n}_{|U^c}\ll \nu,$$
holds true, for all indexes $n$: indeed, if $E_n$ denotes the element of $\mathcal{A}$ such that ${m_n}_{|E_n}\ll\nu$ and ${m_n}_{|E_n^c}\perp\nu$ (see Theorem \ref{troncate}), it 
is enough
 to set $U=\bigcup_n E_n$ (see also {\rm \cite[Theorem 5.2 and following notes]{candeloro}}).

Similarly, in the same situation, there exists a {\em unique} set $V\in \As$ such that

$$m_n^a= {m_n}_{|V},\ \ m_n^c={m_n}_{|V^c},$$
holds true, for all $n's$. So, if the sets $U$ and $V$ in $\As$ are chosen in order to get the decompositions above (including in the sequence $(m_n)$ also the limit $m$), the conclusion of the theorem is a direct consequence of the convergence of the measures ${m_n}_{|U},\ {m_n}_{|U^c} ,\ {m_n}_{|V}, \ {m_n}_{|V^c}$ respectively to ${m}_{|U},\ {m}_{|U^c} ,\ {m}_{|V}, \ {m}_{|V^c}$. \ \ $\Box$

%________________________________________________________________________________________
%_______________________________________________________________________________________-
\subsection{ Convergence for decompositions in the finitely-additive case}
We first establish a convergence result, without assuming particular conditions on the filter $\Fs$.

\medskip

Assume that $(m_n)_n$ is any sequence of uniformly bounded ${\bf X}$-valued finitely additive measures defined on a finitely additive measure space $(\Omega,\Hs,\nu)$. In order to deduce our result, we shall make use of Theorem \ref{stone}.

\begin{definizione}\rm
We say that a sequence $(m_n)_n$ of ${\bf X}$-valued finitely additive measures, defined on a algebra $\Hs$, is {\em uniformly $s$-bounded} if there exists an $(o)$-sequence $(r_k)_k$ in ${\bf X}$ such that, for every disjoint sequence $(H_j)$ from $\Hs$  and for every $k\in \enne$ a corresponding index $j(k)$ can be found, such that
$$\sup_{n\in \mathbb{N}} \sup_{j\geq j(k)}|m_n(H_j)|\leq r_k.$$
\end{definizione}

For example, the sequence  $(m_n)_n$  is uniformly $s$-bounded if there exists an $s$-bounded positive finitely additive measure $m$ such that $|m_n|\leq m$ holds, for every $n$.

\begin{osservazione}{\label{uniformestension}}\rm
We remark here that, in case the finitely additive measures $m_n$ are uniformly $s$-bounded, then it is possible to apply Theorem \ref{stone}, and the density property {\bf \ref{stone}.2} holds in a uniform way, with respect to a unique regulating $(o)$-sequence, see {\rm \cite[Theorem 4.4]{prace}}.
\end{osservazione}

Now we are in position to state the convergence theorem.

\begin{teorema}{\label{main}}
Assume that the ${\bf X}$-valued measures $(m_n)_n$, defined on  $\Hs$, are {\em uniformly $s$-bounded} and pointwise $(o_{\mathcal F})$-convergent to some finitely additive measure $m$. Assume also that $\nu:\Hs\to \erre^+_0$ is a fixed finitely additive  measure.
Then the sequences
$$(m_n^{<})_n,\ (m_n^{\perp})_n, \ (m_n^a)_n,\ \ (m_n^s)_n$$ 
converge in the same way to
$$m^{<},\ m^{\perp}, \ m^a,\ \ m^s$$ 
 respectively, where absolute continuity and singularity are meant w.r.t. $\nu$.
\end{teorema}

\begin{osservazione}\rm
In general, filter-convergence does not imply uniform $s$-boundedness, unless the filter is $\Fs_{cofin}$ and $\Hs$ is a $\sigma$-algebra.  
\end{osservazione}

{\bf Proof of Theorem \ref{main}} \ We first remark that the limit measure $m$ is $s$-bounded too, since the sequence  $(m_n)_n$ is uniformly $s$-bounded.  Using the notation in 
{\bf ($H_{\mbox{\tiny Stone}}$)} 
and Theorem \ref{stone}, and setting
$$\mu_n:=m_n\cdot \psi^{-1},\ \  \mu:=m\cdot \psi^{-1}$$
for each $n$, then $\mu_n$ is defined on $\Sigma$, is $\sigma$-additive in that algebra, (and the same holds for $\mu$), and the sequence $(\mu_n)_n$ is uniformly $s$-bounded. 
Then, thanks to Remark \ref{uniformestension}, 
 the measures $m_n$ can be extended to $\sigma$-additive measures $\widetilde{\mu_n}$ in a uniform way in the $\sigma$-algebra $\sigma(\Sigma)$. Of course, we also have that the sequence $(\mu_n)_n$ is $(o_{\mathcal F})$-convergent to $\mu$ in the algebra $\Sigma$.

Our next step is to prove that the $\sigma$-additive measures $\widetilde{\mu_n}$ satisfy the hypotheses of Theorem \ref{teokyber}. To this aim, let us denote with $\widetilde{\mu}$ the measure extending $\mu$ in $\sigma(\Sigma)$ and let us prove that the sequence  
$\widetilde{\mu_n}$ is $(o_{\mathcal F})$-convergent to $\widetilde{\mu}$.

Let $(p_n)_n$ denote any $(o)$-sequence regulating  pointwise $(o_{\mathcal F})$-convergence: thus for every set $A\in\Sigma$ there exists a sequence $(F_n)_n$ in $\Fs$ such that
$$\sup_{k\in F_n}|\mu_k(A)-\mu(A)|\leq p_n$$
holds true, for all $n$. 

Let us now fix any element $H\in\sigma(\Sigma)$. Thanks to the uniform extendability (see Remark \ref{uniformestension}),
there exist an $(o)$-sequence $(r_j)_j$ in ${\bf X}$ (independent of $H$) and a corresponding sequence $(A_j)_j$ in $\Sigma$, such that
$$\sup_{k\in \enne}|\widetilde{\mu}_k(H)-\mu_k(A_j)|\leq r_j$$
holds for all $j$. (In the last formula, also the measures $\mu$ and $\widetilde{\mu}$ can be included in the $\mu_k's$ and $\widetilde{\mu}_k's$). Now, if we fix arbitrarily any index $j$, for every integer $n\geq j$ there exists $F_n\in \Fs$ such that
$$\sup_{k\in F_n}|\mu_k(A_j)-\mu(A_j)|\leq p_n\leq p_j$$
so, for any $k\in F_n$ we get
$$|\widetilde{\mu}_k(H)-\widetilde{\mu}(H)|\leq |\widetilde{\mu}_k(H)-\mu_k(A_j)|+|\mu_k(A_j)-\mu(A_j)|+|\mu(A_j)-\widetilde{\mu}(H)|\leq 2r_j+p_n\leq 2r_j+p_j.$$
This shows that, choosing the $(o)$-sequence $(q_j:=2r_j+p_j)$ (not depending on $H$), for every index $j$ an element $F\in \Fs$ can be found, such that
$$\sup_{k\in F}|\widetilde{\mu}_k(H)-\widetilde{\mu}(H)|\leq q_j,$$
i.e. the sequence of measures $(\widetilde{\mu}_k)_k$ is $(o_{\mathcal F})$-convergent to $\widetilde{\mu}$. \\
This is clearly sufficient, thanks to Theorem \ref{teokyber}, to deduce also convergence of the decompositions of the measures $m_n$. $\Box$\\

In a similar way we can obtain a convergence result  involving the Yosida-Hewett decompositions.
%________________________________________________________________________________________
%_______________________________________________________________________________________-
\section{Uniform $s$-boundedness}
In order to obtain a result of uniform $s$-boundedness of the sequence $(m_n)_n$, we shall make use of a Schur-type result, proved in \cite{bdpschurfilters}. 
Indeed, in order to prove our Theorem \ref{finale} we shall proceed by contradiction, thus giving rise to a sequence $(H_n)_n$ of pairwise disjoint elements of the $\sigma$-algebra $\mathcal{A}$, along which uniform $s$-boundedness is violated: these sets will be considered as {\em atoms} of a $\sigma$-algebra isomorphic with $\Ps(\enne)$, where the mentioned Schur-type Theorem can be applied in order to get an absurd.
\subsection{A Schur-type theorem for filter convergence in vector lattices}
 
 Here we assume that $\Omega=\enne$, $\Hs=\Ps(\enne)$, and 
 $(m_n)_n$ is a sequence of ${\bf X}$-valued {\em $\sigma$-additive} measures on $\Ps(\enne)$.
 
The following result is a Schur-type theorem for filter convergence, in our setting, see \cite{bdpschurfilters}:  
we just remark once more that, since ${\bf X}$ is super-Dedekind complete and weakly-$\sigma$-distributive, formulations involving  $(o)$-sequences $(p_n)_n$ are equivalent to formulations involving $D$-sequences $(a_{i,j})$.

\begin{teorema}{\label{schur}}{\rm (\cite[Theorem 3.1]{bdpschurfilters})}
 Assume that $\Fs$ is block-respecting and diagonal, and that the measures $m_n$ are uniformly bounded and are $(o_{\mathcal F})$-convergent to 0 with respect to the same regulating $(o)$-sequence $(b_j)_j$.

Then there exists 
an $(o)$-sequence $(a_j)_j$ regulating $\sigma$-additivity of all the measures $(m_n)_n$ and $m$, and, setting $p_j:=2(a_j+2b_j)_j$, an $(o)$-sequence is obtained, with the property that
%an $(o)$-sequence $(p_j)_j$ such that, 
for every integer $j$   there exists $F\in \Fs$ such that for every $n\in F$ the following inequality holds:
$$\sum_{k\in \enne}|m_n(\{k\})|\leq p_j,$$
 i.e.
the sequence $n\mapsto \sum_{k\in \enne}|m_n(\{k\})|$ is $(o_{\mathcal F})$-convergent to $0$.
\end{teorema}
\begin{comment}
\begin{osservazione}\label{notaschur} \rm
The above sequence $(p_j)_j$ can be obtained as follows:
$$p_j=2(a_j+2b_j),$$
where $(a_j)_j$ is an $(o)$-sequence regulating $\sigma$-additivity of all measures $m_n$ (independence of this $(o)$-sequence from $n$ follows from equiboundedness of the measures).
\end{osservazione}
\end{comment}
%________________________________________________________________________________________
%_______________________________________________________________________________________-
\subsection{Filter convergence and uniform $s$-boundedness }
Here we turn back to the general case, and consider finitely additive mesures defined on an arbitrary $\sigma$-algebra $\As$ in a general space $\Omega$.
 In order to deduce a result of uniform $s$-boundedness, the following partial condition will be requested here.

\begin{definizione}{\label{idealsbound}}
Given a sequence  $\{m_j:j\in \enne\}$ of $s$-bounded finitely additive measures on $\Hs$, and a filter $\Fs$ in $\Ps(\enne)$, we say that the measures  $\{m_j:j\in \enne\}$ are {\em ideally uniformly $s$-bounded}  if there exists an $(o)$-sequence $(r_k)_k$ such that, for any family $(H_l)_l$ of pairwise disjoint sets in $\Hs$,
 any integer $k$  and any element $I$ of the dual ideal of $\Fs$, there exists an integer $l(k)$ such that
$$\sup_{j\in I}\sup_{l\geq l(k)}|m_j(H_l)|\leq r_k.$$
\end{definizione}
Of course, if the filter is $\Fs_{cofin}$, the dual ideal consists of all finite sets in $\enne$, hence ideal uniform $s$-boundedness is equivalent to  $s$-boundedness of all the measures $m_j$, provided they are equibounded. \\

Moreover, ideal  uniform $s$-boundedness of an equibounded sequence of measures $(m_n)_n$ 
in a countably generated filter (see Definition \ref{4.3})
can be easily obtained if they are each $s$-bounded with respect to the same regulator and if uniform $s$-boundedness holds separately for the subsequences $(m_j)_{j \in A_k}$.\\
%_________________________________________

Now,  we shall prove that, if the filter $\Fs$ enjoys particular properties, ideal uniform $s$-boundedness is sufficient to ensure uniform $s$-boundedness of the  sequence $(m_n)_n$. As usual, we shall assume that $\Fs$ is a fixed  filter in $\enne$, so ideal uniform $s$-boundedness is related to the dual ideal $\mathcal{I}_{\mathcal{F}}$.
The theorem is as follows.
\begin{teorema}{\label{finale}}
Assume that $(m_n)_n$ is an equibounded sequence of ideally uniformly $s$-bounded finitely additive measures, defined on a $\sigma$-algebra $\As$, and taking values in ${\bf X}$. If the filter $\Fs$ is block-respecting and diagonal, and if the measures $m_n$ are $(o_{\mathcal F})$-convergent {\rm(}with respect to the same regulating $(o)$-sequence{\rm)} to an $s$-bounded finitely additive measure $m$, then the measures are uniformly $s$-bounded.  
\end{teorema}

{\bf Proof.}\  Let $(q_j)_j$ be any $(o)$-sequence in ${\bf X}$, regulating ideal uniform $s$-bound\-ed\-ness of the measures $m_n$ (including also the limit measure $m$). Of course, $(q_j)_j$ regulates also simple $s$-boundedness of each measure $m_n$. Moreover, let us denote by $(b_j)_j$ an $(o)$-sequence regulating $(o_{\mathcal F})$-convergence. Without loss of generality, we can and do assume that the same $(o)$-sequences regulate ideal uniform $s$-boundedness and convergence of the measures $(m_n-m)$, so that we can replace $m$ with 0.
We claim that the $(o)$-sequence $(r_j)_j$, defined as
$$r_j=2(q_j+2b_j)+q_j$$
regulates uniform $s$-boundedness of the measures $m_n$.

Otherwise, by a usual procedure, one could find an integer $j_0$, a  pairwise disjoint family $(H_k)_k$ in $\As$, and a corresponding sequence $(n_k)_k $ such that
\begin{eqnarray}{\label{contra}}
|m_{n_k}(H_k)|\not\leq r_{j_0}\end{eqnarray}
for all $k$.

But the family $\{H_k:k\in \enne\}$ has a subsequence $\{H_{k_l}:l\in \enne\}$ along which all the measures $m_n$ are $\sigma$-additive, with $(q_j)_j$ as regulating $(o)$-sequence (see Remark \ref{processodiagonale}). Thus, if we consider the sets $H_{k_l}$ 
as atoms, the mappings $m_n$ can be viewed as $\sigma$-additive measures defined on the $\sigma$-algebra $\Hs$ generated by these atoms (isomorphic with $\Ps(\enne)$). 
Of course, these measures restricted to $\Hs$ satisfy the hypotheses of Theorem \ref{schur}, and therefore their total variations $(o_{\mathcal F})$-converge uniformly to 0, with respect to the $(o)$-sequence $2(q_j+2b_j)_j$.
This means that, even for the exceptional integer $j_0$, an element $F_0\in \Fs$ exists, such that the set 
$\{n\in \enne: \sup_{l}|m_n(H_{k_l})|\leq 2(q_{j_0}+2b_{j_0})\}$ contains $F_0$.

Now, the set $F_0^c$ belongs to the dual ideal, so, by ideal uniform $s$-boundedness, there exists $l_0$ such that
$$\sup_{l\geq l_0}\sup_{n\in F_0^c}|m_n(H_{k_l})|\leq q_{j_0}.$$
So, for every integer $n$ and every $l\geq l_0$ we have $|m_n(H_{k_l})|\leq r_{j_0}$,
and in particular 
$$|m_{n_{k_l}}(H_{k_l})|\leq r_{j_0},$$
thus contradicting (\ref{contra}). \ \ $\Box$

\medskip

Putting together Theorem \ref{finale} and Theorem \ref{main}, we obtain the following result on convergence of decompositions. 

\begin{corollario}{\label{corofinale}}
Assume that $(m_n)_n$ is an equibounded sequence of ideally uniformly $s$-bounded finitely additive measures, defined on a $\sigma$-algebra $\As$, and taking values in ${\bf X}$. Assume that the filter $\Fs$ is block-respecting and diagonal, and that the measures $m_n$ are $(o_{\mathcal F})$-convergent to an $s$-bounded finitely additive measure $m$. 
Moreover, let $\nu:\As\to \erre^+_0$ be any positive finitely additive measure. 
Then the sequences
$$(m_n^{<})_n,\ (m_n^{\perp})_n, \ (m_n^a)_n,\ \ (m_n^s)_n$$ converge in the same way to
$$m^{<},\ m^{\perp}, \ m^a,\ \ m^s$$ 
 respectively, where absolute continuity and singularity are meant w.r.t. $\nu$.
\\
Finally, also the Yosida-Hewett decompositions of the measures $m_n$ $(o_{\mathcal F})$-converge to the Yosida-Hewett decomposition of $m$.
\end{corollario}

Next, we shall prove an alternative version of  Theorem \ref{finale},
relaxing some conditions on the domain of the measures but strenghtening conditions on the filter.
We shall assume that the filter $\Fs$ satisfies the following condition.
\begin{definizione}\label{fortediagonale}
A filter $\Fs$ in $\Ps(\enne)$ is said to be {\em super-diagonal} if every stationary set $I$ has a stationary subset $J$ such that all infinite subsets of $J$ are stationary. 
\end{definizione}
It is not difficult to see that a super-diagonal filter is always diagonal and block-respecting. 

However, every countably generated filter is super-diagonal, as proved in Proposition \ref{superdiagonal}.
 
On the other hand,
we only require that the measures $m_n$ are uniformly bounded, finitely additive and $s$-bounded. Moreover, we shall relax the condition of convergence, just assuming 
$(o_{\mathcal F})$-convergence to 0 in some special algebra, with respect to the same regulating $(o)$-sequence. We begin with some definitions.
\begin{definizione}\rm
Let $\Hs$ be an algebra of subsets of an abstract set $\Omega$. We say that $\Hs$ enjoys the {\em property (SCP)} if, for every sequence $(H_k)_k$ of pairwise elements from $\Hs$, there exists a subsequence $(H_{k_r})_r$ whose union belongs to $\Hs$.
\end{definizione}

We remark that $(SCP)$ is one of the conditions on an algebra that  ensure the validity of the Vitali-Hahn-Saks and related theorems in the finitely additive case. We  refer to \cite{haydon,freniche,mimmo} for alternative conditions and further results.

\begin{definizione}\rm 
Given an algebra $\Hs$  as above, and a sequence of $s$-bounded finitely additive measures $m_n$ on it, we say that the $m_n's$ are {\em $\Fs$-uniformly $s$-bounded} if there exists an $(o)$-sequence $(r_p)_p$ such that, for every sequence $(H_k)_k$ of pairwise disjoint elements from $\Hs$ and every integer $p$ there exist an integer $k(p)$ and an element $F(p)\in \Fs$ such that $|m_n(H_k)|\leq r_p$ for all $k\geq k(p)$ and all $n\in F(p)$. 
\end{definizione}

The result reads as follows.
\begin{teorema}\label{nuovoschur}
Let us assume that a sequence $(m_n)_n$ of uniformly bounded, $s$-bounded finitely additive measures is defined on an algebra $\Hs$  of subsets of $\Omega$, enjoying the property {\em (SCP)}. Moreover, assume that the filter $\Fs$ is super-diagonal and that
there exists an $(o)$-sequence $(b_j)_j$ such that for every set $A\in \Hs$ and every integer $j$ there exists $F\in \Fs$ such that 
$|m_n(A)|\leq b_j$
for every $n\in F$.
Then the measures $(m_n)_n$ are $\Fs$-uniformly $s$-bounded.
\end{teorema}
{\bf Proof.}\ Let $(q_j)_j$ be any $(o)$-sequence regulating $s$-boundedness of all the measures $v(m_n)$, and let
 $u:=\sup_nv(m_n)(\Omega)$. Thanks to a result from
{\rm \cite[Lemma 2.7]{panoramica}},  there exists an increasing
 map $w:\enne\to \enne$ such that the following setting defines two $(o)$-sequences:
$$Q_N:=u\wedge \sum_{j=N}^{\infty}q_{w(j)},\ \ B_N:=u\wedge \sum_{j=N}^{\infty}b_{w(j)}.$$
Now, let $r_p:=2(b_p+B_p+Q_p)$ for all $p\in\enne$, and let us
 prove that for this $(o)$-sequence the following claim is true:

\begin{description}
\item[\ref{nuovoschur}.0] 
for every disjoint sequence $(H_j)_j$ in $\Hs$ and every integer $p$ the set
$$\{j\in \enne: \bigvee_{k\in \enne}|m_j(H_k)|\leq r_p\}$$
belongs to $\Fs$. 
\end{description}
This is actually a stronger condition than $\Fs$-uniform $s$-boundedness with respect to $\Hs$. 

We shall proceed by contradiction: if  
{\bf \ref{nuovoschur}.0}
is false, then there exist an integer $p_0$ and a disjoint sequence $(H_j)_j$ in $\Hs$ such that the set
$$I:=\{j\in \enne: \bigvee_{k\in \enne}|m_j(H_k)|\not\leq r_{p_0}\}$$
is stationary. Since $\Fs$ is diagonal, then there exists a stationary set $J\subset I$ such that, for every integer $p$ and every index $k$ an integer $j(p,k)$ exists, for which
$$|m_j(H_k)|\leq b_p$$
holds, as soon as $j\geq j(p,k)$, $j\in J$ (see {\rm \cite[Lemma 2.2]{bdpschurfilters}}). Moreover, thanks to super-diagonality, we can assume that every infinite subset of $J$ is stationary.

Now, for every integer $j\in J$ an element $H_{k_j}$ exists, such that  
$$|m_j(H_{k_j})|\not\leq r_{p_0}.$$
From now on, we shall denote by $H'_j$ the set $H_{k_j}$, $j=1,2...$. So, we can write
\begin{eqnarray}\label{assurdo}
|m_j(H'_j)|\not\leq r_{p_0},\end{eqnarray}
for all $j\in J$.

Let us denote by $j_1$ the first element of $J$. Thanks to the $(SCP)$ property and to $s$-boundedness of $m_{j_1}$, by a
 standard procedure it is possible to find an infinite subset $P_1\subset J$ such that
$$v(m_{j_1})(E_1)\leq q_{w(p_0+1)},$$
where $E_1=\bigcup_{j\in P_1}H'_j$.
Let now $p_1$ denote the minimum of $P_1$. Thanks to the mentioned consequence of diagonality, we see that there exists an integer $N_1>p_1$, $N_1\in J$, such that 
$$|m_i(H'_{p_1})|\leq b_{w(p_0+1)},$$
for all $i\in J,\ i\geq N_1$.

Again, thanks to the $(SCP)$ property and to $s$-boundedness of
 $m_{p_1}$, it is possible to find an infinite subset $P_2\subset P_1$ such that
$$v(m_{p_1})(E_2)\leq q_{w(p_0+2)},$$
where $E_2=\bigcup_{j\in P_2}H'_j$.
Without loss of generality, we shall assume that the minimum element of $P_2$, which we denote by $p_2$, is strictly larger than $N_1+p_1.$

As above, thanks to the mentioned consequence of diagonality, we see that there exists an integer $N_2>p_2$, $N_2\in J$, such that 
$$|m_i(H'_{p_2})|\leq b_{w(p_0+2)},$$
for all $i\in J,\ i\geq N_2$.

Proceeding in this way, we obtain a decreasing sequence $(P_s)_s$ of infinite subsets of $J$, a corresponding decreasing sequence  $(E_s)_s$ in $\Hs$, $E_s=\bigcup_{i\in P_s}H'_i$, and two increasing sequences $(N_s)_s,\ (p_s)_s$ in $J$, such that, for all $s$:
\begin{description}
\item[\ref{nuovoschur}.1]  $N_s>p_s,\ \ p_{s+1}>N_s+p_s$;

\item[\ref{nuovoschur}.2]  $H'_i\subset  E_s$ whenever $i\in P_s$; \ $H'_i\cap E_{s+1}=\emptyset$ whenever $i\leq p_s$;

\item[\ref{nuovoschur}.3] $\bigvee_{i=1}^{p_s}v(m_i)
(E_{s+1})\leq u\wedge \sum_{i=1}^{p_s}q_{w(p_0+i+1)}\leq Q_{p_0}$;

\item[\ref{nuovoschur}.4]  $|m_i(H'_{p_s})|\leq u\wedge b_{w(p_0+s)}$, whenever $i\geq N_s,\ i\in J.$
\end{description}
Now, the sequence $(H'_{p_s})_s$ has a subsequence $(H'_{{p_{s_l}}})_l$ whose union $H$ belongs to $\Hs$.
Since the sequence $(m_n)_n$ is $(o_{\mathcal F})$-convergent to 0, there exists an element $F\in \Fs$ such that $|m_i(H)|\leq b_{p_0}$ for all $i\in F$. 
Since set $J':=\{p_{s_l}:l\in \enne\}$ is an infinite subset of $J$, it is stationary, hence has an element in common with $F$: let $p_{s_{l'}}$ be such element, and denote by simplicity with $m^*$ the measure $m_{p_{s_{l'}}}$ and by $H^*$ the set $H_{p_{s_{l'}}}$ . From (\ref{assurdo}) we should have that $|m^*(H^*)|\not\leq r_{p_0}$.
But we see that:
$$m^*(H^*)=m^*(H)-m^*(\bigcup_{i<l'}H'_{p_{s_i}})-m^*(\bigcup_{i> l'}H'_{p_{s_i}}). $$
Since the set $\bigcup_{i> l'}H'_{p_{s_i}}$
 is contained in $E_{s_{l'}+1}$, we obtain 
$$v(m^*)(\bigcup_{i> l'}H'_{p_{s_i}})\leq Q_{p_0} $$
thanks to  {\bf \ref{nuovoschur}.3}.  
Moreover, since $p_{s_{l'}}>N_{s_{l'}}$, from 
 {\bf \ref{nuovoschur}.4} we deduce 
$$|m^*(\bigcup_{i<l'}H'_{p_{s_i}})|\leq u\wedge \sum_{i<l'}b_{w(p_0+s_i)}\leq B_{p_0}.$$
Finally we have $|m^*(H)|\leq b_{p_0}$ since $p_{s_{l'}}\in F$.
Thus we obtain:
$$|m^*(H^*)|\leq b_{p_0}+Q_{p_0}+B_{p_0}\leq r_{p_0},$$
but this is in contrast with formula (\ref{assurdo}). \ \ $\Box$

Finally, from Theorem \ref{nuovoschur} and Theorem \ref{main}, we obtain the following result on convergence of decompositions. 

\begin{corollario}\label{nuovocorfinale}
Assume that $(m_n)_n$ is an equibounded sequence of ideally uniformly $s$-bounded finitely additive measures, defined on an algebra $\As$ enjoying $(SCP)$, and taking values in ${\bf X}$. Assume that the filter $\Fs$ is  super-diagonal, and that the measures $m_n$ are $(o_{\mathcal F})$-convergent to an $s$-bounded finitely additive measure $m$. 
Moreover, let $\nu:\As\to \erre^+_0$ be any positive finitely additive measure. 
Then the sequences
$$(m_n^{<})_n,\ (m_n^{\perp})_n, \ (m_n^a)_n,\ \ (m_n^s)_n$$ converge in the same way to
$$m^{<},\ m^{\perp}, \ m^a,\ \ m^s$$ 
 respectively, where absolute continuity and singularity are meant w.r.t. $\nu$.
\\
Finally, also the Yosida-Hewett decompositions of the measures $m_n$ $(o_{\mathcal F})$-converge to the Yosida-Hewett decomposition of $m$.
\end{corollario}
{\bf Proof.}\ Thanks to Theorem \ref{nuovoschur}, the measures $m_n$ are $\Fs$-uniformly $s$-bounded. From this, and from ideal uniform $s$-boundedness, it follows easily uniform $s$-boundedness. Then, from Theorem \ref{main} we get the assertion. \ \ $\Box$

\begin{comment}
\bibitem{AV1} Anna Avallone - Paolo Vitolo, \textit{DECOMPOSITION AND CONTROL
THEOREMS IN EFFECT ALGEBRAS}, Scientiae Mathematicae Japonicae Online,
 Vol. 8, (2003), 1-14.
\bibitem{ABVW1} A. Avallone -G. Barbieri - P. Vitolo - H. Weber,
 \textit{Decomposition of effect algebras and the Hammer-Sobczyk theorem},
 Algebra Universalis {\bf 60},   (2009) 1-18. DOI 10.1007/s00012-008-2083-z.
\bibitem{AV2} A. Avallone - P. Vitolo, \textit{Lyapunov decomposition of
 measures on effect algebras}, Sci. Math. Japonicae {\bf 69}, No. 1,  (2009) 79-87.
\bibitem{W1} H. Weber, \textit{Boolean algebras of lattice uniformities and
 decompositions of modular functions}, Ricerche di Matematica
 {\bf 58},  (2009) 15-32 . DOI 10.1007/s11587-009-0043-3
\bibitem{ABV1} A. Avallone - G. Barbieri - P. Vitolo, \textit{Central
 elements in pseudo-D-lattices and Hahn decomposition theorem}
 Bollettino UMI, III, {\bf 9}, no. 3,  (2010) 447-470.
\end{comment}

\end{document}